\documentclass[11pt]{amsart}
\usepackage{amssymb}
\usepackage{amsmath}
\usepackage{amscd}
\usepackage[matrix,arrow]{xy}
\usepackage{url}
\usepackage{color}
\usepackage{xtab}
\newcommand{\Q}{{\mathbb Q}}
\newcommand{\Z}{{\mathbb Z}}
\newcommand{\N}{{\mathbb N}}

\newcommand{\R}{{\mathbb R}}

\newcommand{\op}{\operatorname}

\def\ellog{{\mathfrak{l}}}

\begin {document}
\bibliographystyle{plain}
\newtheorem{thm}{Theorem}

\newtheorem*{conj}{Conjecture}

\theoremstyle{definition}

\theoremstyle{remark}

\newtheorem*{rem}{Remark}

\title[]{On the Diophantine equation $\binom{n}{k}=\binom{m}{l}+d$}

\author{H. R. Gallegos-Ruiz}

\address{Homero R. Gallegos-Ruiz\\
	Unidad Acad\'emica de Matem\'aticas\\
	Universidad Aut\'onoma de Zacatecas\\
	Calzada Solidaridad y Paseo de la Bufa\\
	Zacatecas, Zacatecas, CP 98000 \\Mexico
}
\email{hgallegos@uaz.edu.mx}


\author{N. Katsipis}
\address{Nikolaos Katsipis\\ Department of Mathematics \& Applied Mathematics\\ University of Crete\\ GR-70013, Heraklion, Crete\\ Greece}
\email{katsipis@gmail.com}

\author{Sz.\ Tengely}
\address{Szabolcs Tengely\\
Institute of Mathematics\\
University of Debrecen\\
P.O.Box 12\\
4010 Debrecen\\
Hungary}

\email{tengely@science.unideb.hu}

\author{M. Ulas}
\address{
Maciej Ulas\\
Jagiellonian University\\
Faculty of Mathematics and Computer Science \\
Institute of Mathematics\\
\L{}ojasiewicza 6\\
30-348 Krak\'ow\\
Poland}

\address{and}

\address{Institute of Mathematics of the Polish Academy of Sciences \\
\'{S}wi\c{e}tego Tomasza 30\\
 31-014 Krak\'{o}w, Poland}
\email{Maciej.Ulas@im.uj.edu.pl}

\date{\today}
\thanks{}

\keywords{binomial coefficient, Diophantine equation, elliptic curve, genus two curve, integer points}
\subjclass[2000]{Primary 11G30, Secondary 11J8}

\begin{abstract}
By finding all integral points on certain elliptic and hyperelliptic curves we completely solve the Diophantine equation $\binom{n}{k}=\binom{m}{l}+d$ for $-3\leq d\leq 3$ and $(k,l)\in\{(2,3),\; (2,4),\;(2,5),\; (2,6),\; (2,8),\; (3,4),\; (3,6),\; (4,6), \; (4,8)\}.$ Moreover, we present some other observations of computational and theoretical nature concerning the title equation. 
\end{abstract}
\maketitle

\section{Introduction}
There are many nice results related to the equation
\begin{equation}\label{binom}
\binom{n}{k}=\binom{m}{l},
\end{equation}
in unknowns $k$, $l$, $m$, $n.$
This is usually considered with the restrictions
$2\leq k\leq n/2, 2\leq l\leq m/2$ and $k<l.$
The only known  solutions (with the above mentioned restrictions)
are the following
\begin{eqnarray*}
	&& \binom{16}{2}=\binom{10}{3},\quad \binom{56}{2}=\binom{22}{3},\quad \binom{120}{2}=\binom{36}{3},\\
	&& \binom{21}{2}=\binom{10}{4},\quad \binom{153}{2}=\binom{19}{5},\quad \binom{78}{2}=\binom{15}{5}=\binom{14}{6},\\
	&& \binom{221}{2}=\binom{17}{8},\quad \binom{F_{2i+2}F_{2i+3}}{F_{2i}F_{2i+3}}=\binom{F_{2i+2}F_{2i+3}-1}{F_{2i}F_{2i+3}+1} \mbox{ for } i=1,2,\ldots,
\end{eqnarray*}
where $F_n$ is the $n$th Fibonacci number. The infinite family of solutions involving Fibonacci numbers was found by Lind \cite{Lind} and Singmaster \cite{Singmaster}.

Equation \eqref{binom} has been completely solved for pairs
\[
(k,l)=(2,3),\; (2,4),\; (2,6),\; (2,8),\; (3,4),\; (3,6),\; (4,6), \; (4,8).
\]
In cases of these pairs one can easily reduce the equation to the determination of solutions of a number
of Thue equations or elliptic Diophantine equations.
In 1966, Avanesov \cite{Avanesov} found all integral solutions of
equation \eqref{binom} with $(k,l)=(2,3).$
De~Weger \cite{deWeger} and independently
Pint\'er \cite{Pinter} provided all the solutions of the equation with $(k,l)=(2,4).$
The case $(k,l)=(3,4)$ reduces to the equation
$Y(Y+1)=X(X+1)(X+2)$ which was solved by Mordell \cite{Mordell}.
The remaining pairs $(2,6),(2,8),(3,6),(4,6),(4,8)$ were
handled by Stroeker and de Weger \cite{StroekerWeger}, using linear forms in elliptic logarithms.
The case with $(k,l)=(2,5)$ was completely solved by Bugeaud, Mignotte, Siksek, Stoll and Tengely \cite{BMSST}, the integral solutions are as follows
\begin{gather*}
(n,m)=(0, 0), \;
(0, 1),\;
(1, 0),\;
(1, 1),\;
(2, 0),\;
(2, 1),\;
(3, 0),\;
(3, 1),\;
(4, 0),\;
(4, 1),\;
(5, -1),\\
(5, 2),\;
(6, -3),\;
(6, 4),\;
(7, -6),\;
(7, 7),\;
(15, -77),\;
(15, 78),\;
(19, -152),\;
(19, 153).
\end{gather*}
In a recent paper Blokhuis, Brouwer and de Weger \cite{BBW} determined all non-trivial
solutions with $\binom{n}{k}\leq 10^{60}$ or $n\leq 10^6.$ General finiteness results are also known.
In 1988, Kiss \cite{Kiss} proved that if $k=2$ and $l$ is a given odd prime, then the equation has only finitely
many positive integral solutions. Using Baker's method, Brindza \cite{Brindza2}
showed that equation \eqref{binom} with $k=2$ and $l\geq 3$ has only finitely many positive integral solutions.

In case of the more general equation
\begin{equation}\label{binomd}
\binom{n}{k}=\binom{m}{l}+d
\end{equation}
Blokhuis, Brouwer and de Weger \cite{BBW} determined all non-trivial
solutions with $d=1$ and $(k,l),(l,k)=(2,3),(2,4),(2,6),(3,4),(4,6),(4,8)$ and $(k,l)=(2,8).$ They provided a complete list of solutions for the above cases and if $\binom{n}{k}\leq 10^{30}.$
{\begin{center}
$\begin{array}{ll}\small
\begin{tabular}{|l|l|l|l|} \hline
	$n$ & $k$ & $m$ & $l$ \\ \hline \hline
	$11$ & $2$ & $8$ & $3$ \\ \hline
	$60$ & $2$ & $23$ & $3$ \\ \hline
	$160403633$ & $2$ & $425779$ & $3$ \\ \hline
	$6$ & $3$ & $7$ & $2$ \\ \hline
	$7$ & $3$ & $9$ & $2$ \\ \hline
	$16$ & $3$ & $34$ & $2$ \\ \hline
	$27$ & $3$ & $77$ & $2$ \\ \hline
	$29$ & $3$ & $86$ & $2$ \\ \hline
	$34$ & $3$ & $21$ & $4$ \\ \hline
\end{tabular}
&\small
\begin{tabular}{|l|l|l|l|} \hline
$n$ & $k$ & $m$ & $l$ \\ \hline \hline
$19630$ & $3$ & $1587767$ & $2$ \\ \hline
$12$ & $4$ & $32$ & $2$ \\ \hline
$93$ & $4$ & $2417$ & $2$ \\ \hline
$10$ & $5$ & $23$ & $2$ \\ \hline
$22$ & $5$ & $230$ & $2$ \\ \hline
$62$ & $5$ & $3598$ & $2$ \\ \hline
$135$ & $5$ & $26333$ & $2$ \\ \hline
$139$ & $5$ & $28358$ & $2$ \\ \hline
$28$ & $11$ & $6554$ & $2$ \\ \hline
\end{tabular}
\end{array}$
\end{center}}
\begin{center}
Table 1. Known solutions of the Diophantine equation $\binom{n}{k}=\binom{m}{l}$. 
\end{center}
If $d$ is not fixed they also obtained some interesting infinite families, an example is given by
$$
\binom{12x^2-12x+3}{3}+\binom{x}{2}=\binom{24x^3-36x^2+15x-1}{2}.
$$

In 2019, Katsipis \cite{Katsipis} completely resolved the case with $(k,l)=(8,2)$ and he also determined the integral solutions if $(k,l),(l,k)=(3,6)$ and $d=1.$

The aim of this paper is to extend results mentioned above and offer some general observations and computational results.

\section{Main results}
	We start our discussion with some numerical observations. More precisely, we observed that for certain pairs $(k, l)$ and an integer $d$, the congruence
\begin{equation}\label{cong}
\binom{n}{k}\equiv \binom{m}{l}+d\pmod{p},
\end{equation}
with suitable chosen prime number $p>\op{max}\{k,l\}$, has no solutions. This immediately implies unsolvability in integers of the related Diophantine equation.

\begin{thm}\label{kl24}
	If $(k,l)=(2,4), d\in\Z$ and $3$ is a quadratic non-residue modulo $p>4$, 
	where the $p$-adic valuation of $12d+1$ is odd,
	then congruence {\rm (\ref{cong})} has no solutions. In particular, equation {\rm (\ref{binomd})} has no solutions in integers.
\end{thm}

\begin{rem}
Based on the previous theorem we may provide some explicit results, for example if $d\equiv u\pmod{75},$ where $u\in\{7, 12, 17, 22, 32, 37, 42, 47, 57, 62, 67, 72\},$ then equation (\ref{binomd}) has no solutions in integers with $(k,l)=(2,4).$
\end{rem}

By using elementary number theory we compute all integral solutions of equation \eqref{binomd} for some values of $k$ and $d$ with $l=k$ and $d\neq 0.$
We note that the case $k=2$ is in some sense trivial. Indeed, in this case the solvability of equation (\ref{keql}) is equivalent to the existence of integers $u,v$ such that $u^2-v^2=8d$ and 
$u\equiv v\equiv 1\pmod{2}$.
Equivalently, we need to determine integers $d_1,d_2$ with $d_{1}\leq d_{2}$ and $8d=d_{1}d_{2}$
satisfying the conditions
$$
d_{1}+d_{2}\equiv 2\pmod{4},\quad d_{2}-d_{1}\equiv 2\pmod{4}.
$$
Thus, if $d$ is odd, one can take $d_{1}=4z_{1}, d_{2}=2z_{2}$, where $d=z_{1}z_{2}$, i.e., the number of solutions of our equation is at least $\sigma_{0}(d)$, where $\sigma_{0}(n)=\sum_{k|n}1$. If $d$ is even one possible choice is $d_{1}=2, d_{2}=4d$.

\begin{thm}\label{kl_eq}
All integral solutions $(n,m)$ of equation \eqref{binomd} with $l=k, k\in\{3,4,5\}$ and $d\neq 0, d\in\{1,2,\ldots,20\}$ are as follows
\begin{center}
\begin{tabular}{|l|l|} \hline
	$(k,d,solutions)$ & $(k,d,solutions)$ \\ \hline \hline
	$\left(3, 3, \left[\left(4, 3\right)\right]\right)$ & $\left(4, 4, \left[\left(5, 4\right)\right]\right)$ \\ \hline
	 $\left(3, 6, \left[\left(5, 4\right)\right]\right)$ & $\left(4, 10, \left[\left(6, 5\right)\right]\right)$\\ \hline
	$\left(3, 9, \left[\left(5, 3\right)\right]\right)$ & $\left(4, 14, \left[\left(6, 4\right)\right]\right)$\\ \hline
	$\left(3, 10, \left[\left(6, 5\right)\right]\right)$ & $\left(4, 20, \left[\left(7, 6\right)\right]\right)$ \\ \hline
	$\left(3, 15, \left[\left(7, 6\right)\right]\right)$ & $\left(5, 5, \left[\left(6, 5\right)\right]\right)$ \\ \hline
	$\left(3, 16, \left[\left(6, 4\right)\right]\right)$ & $\left(5, 15, \left[\left(7, 6\right)\right]\right)$\\ \hline
	$\left(3, 19, \left[\left(6, 3\right)\right]\right)$ & $\left(5, 20, \left[\left(7, 5\right)\right]\right)$\\ \hline

\end{tabular}
\end{center}
\end{thm}
In the next result we deal with the cases that can be reduced to elliptic curves.
\begin{thm}\label{elliptic}
All integral solutions $(m,n)$ of equation \eqref{binomd} with $d\in\{-3,\ldots,3\}$ and $n\geq k, m\geq l$ are as follows.

\begin{center}\small
\begin{tabular}{|c|l|} \hline
	$d$ & $(k,l)=(2,3)$\\ \hline \hline
	$3$ & $\left[\left(75, 368\right), \left(77, 383\right), \left(421726, 158118758\right)\right]$ \\ \hline
	$2$ & $\left[\left(3, 3\right), \left(4, 4\right), \left(104, 604\right)\right]$ \\ \hline
	$1$ & $\left[\left(6, 7\right), \left(7, 9\right), \left(16, 34\right), \left(27, 77\right), \left(29, 86\right), \left(260, 2407\right), \left(665, 9879\right), \left(19630, 1587767\right)\right]$ \\ \hline
	$0$ & $\left[\left(3, 2\right), \left(5, 5\right), \left(10, 16\right), \left(22, 56\right), \left(36, 120\right)\right]$ \\ \hline
	$-1$ & $\left[\left(4, 3\right), \left(8, 11\right), \left(23, 60\right), \left(425779, 160403633\right)\right]$ \\ \hline
	$-2$ & $\left[\right]$ \\ \hline
	$-3$ & $\left[\left(4, 2\right)\right]$ \\ \hline
\end{tabular}
\end{center}

$\begin{array}{ll}\small

\begin{tabular}{|c|l|} \hline
	$d$ &  $(k,l)=(2,4)$\\ \hline \hline
	$3$ & $\left[\right]$ \\ \hline
	$2$ & $\left[(4,3)\right]$ \\ \hline
	$1$ & $\left[(5,4),(7,9),(12,32),(93,2417)\right]$ \\ \hline
	$0$ & $\left[(4,2),(6,6),(10,21)\right]$ \\ \hline
	$-1$ & $\left[\right]$ \\ \hline
	$-2$ & $\left[\right]$ \\ \hline
	$-3$ & $\left[\right]$ \\ \hline
\end{tabular}
&\small
\begin{tabular}{|c|l|} \hline
	$d$ &  $(k,l)=(2,6)$\\ \hline \hline
	$3$ & $\left[(7,5),(11,31),(50,5638)\right]$ \\ \hline
	$2$ & $\left[(6,3)\right]$ \\ \hline
	$1$ & $\left[\right]$ \\ \hline
	$0$ & $\left[(6,2),(8,8),(10,21),(14,78)\right]$ \\ \hline
	$-1$ & $\left[(7,4)\right]$ \\ \hline
	$-2$ & $\left[\right]$ \\ \hline
	$-3$ & $\left[\right]$ \\ \hline
\end{tabular}
\end{array}$

\vspace*{10pt}

$\begin{array}{ll}\small
\begin{tabular}{|c|l|} \hline
$d$ &  $(k,l)=(2,8)$\\ \hline \hline
$3$ & $\left[\right]$ \\ \hline
$2$ & $\left[(8,3)\right]$  \\ \hline
$1$ & $\left[(5,9),(32,12)\right]$ \\ \hline
$0$ & $\left[(8,2),(10,10),(14,78),(17,221)\right]$ \\ \hline
$-1$ & $\left[\right]$ \\ \hline
$-2$ & $\left[\right]$ \\ \hline
$-3$ & $\left[(9,4)\right]$ \\ \hline
\end{tabular}
&\small
\begin{tabular}{|c|l|} \hline
$d$ &  $(k,l)=(3,4)$ \\ \hline \hline
$3$ & $\left[(4,4)\right]$ \\ \hline
$2$ & $\left[\right]$ \\ \hline
$1$ & $\left[\right]$ \\ \hline
$0$ & $\left[(4,3),(7,7)\right]$ \\ \hline
$-1$ & $\left[(5,4),(21,34)\right]$ \\ \hline
$-2$ & $\left[\right]$ \\ \hline
$-3$ & $\left[\right]$ \\ \hline
\end{tabular}
\end{array}$

\vspace*{10pt}

{\small$\begin{array}{lll}
	\begin{tabular}{|c|l|} \hline
$d$ &  $(k,l)=(3,6)$\\ \hline \hline
$3$ &  $\left[(6,4), (7,5)\right]$\\ \hline
$2$ & $\left[\right]$ \\ \hline
$1$ & $\left[\right]$ \\ \hline
$0$ &  $\left[(6,3),(9,9)\right]$ \\ \hline
$-1$ &  $\left[\right]$ \\ \hline
$-2$ & $\left[\right]$ \\ \hline
$-3$ & $\left[(7,4)\right]$ \\ \hline
\end{tabular}
&
\begin{tabular}{|c|l|} \hline
$d$ &  $(k,l)=(4,6)$\\ \hline \hline
$3$ & $\left[\right]$ \\ \hline
$2$ & $\left[\right]$ \\ \hline
$1$ & $\left[\right]$ \\ \hline
$0$ & $\left[(6,4),(10,10)\right]$ \\ \hline
$-1$ & $\left[\right]$ \\ \hline
$-2$ & $\left[(7,5)\right]$ \\ \hline
$-3$ & $\left[\right]$ \\ \hline
\end{tabular}
&
\begin{tabular}{|c|l|} \hline
$d$ &  $(k,l)=(4,8)$\\ \hline \hline
$3$ & $\left[\right]$ \\ \hline
$2$ & $\left[\right]$ \\ \hline
$1$ & $\left[\right]$ \\ \hline
$0$ & $\left[(8,4),(12,12)\right]$ \\ \hline
$-1$ & $\left[\right]$ \\ \hline
$-2$ & $\left[\right]$ \\ \hline
$-3$ & $\left[\right]$ \\ \hline
\end{tabular}
\end{array}$}

\end{thm}
Among the solutions given by Blokhuis, Brouwer and de Weger \cite{BBW} there are some with $(k,l)=(2,5)$ e.g.:
$$
\binom{10}{5}+1=\binom{23}{2},\quad \binom{22}{5}+1=\binom{230}{2},\quad \binom{62}{5}+1=\binom{3598}{2}
$$
in these cases the problem can be reduced to genus 2 curves.
\begin{thm}\label{hyperell}
All integral solutions $(n,m)$ of equation \eqref{binomd} with $d\in\{-3,\ldots,3\},k=2,l=5$ are as follows.
\begin{center}\small
	\begin{tabular}{|c|l|} \hline
		$d$ & solutions \\ \hline \hline
		$-3$ & $\left[(3,6)\right]$ \\ \hline
		$-2$ & $\left[\right]$ \\ \hline
		$-1$ & $\left[(11,8)\right]$ \\ \hline
		$0$ & $\left[(2,5),(4,6),(7,7),(78,15),(153,19)\right]$ \\ \hline
		$1$ & $\left[(23,10),(230,22),(3598,62),(26333,135),(28358,139)\right]$ \\ \hline
		$2$ & $\left[(3,5)\right]$ \\ \hline
		$3$ & $\left[(31,11),(94,16),(346888,375),(356263,379)\right]$ \\ \hline
	\end{tabular}
\end{center}
\end{thm}

Let $k\in\N$ be odd. In the following theorem we consider the Diophantine equation
\begin{equation}\label{k22eq}
\binom{f_{1}(x)}{k}+\binom{x}{2}=\binom{f_{2}(x)}{2}
\end{equation}
in polynomials $f_{1}, f_{2}\in\Q[x]$ satisfying the condition $\op{deg}f_{1}=2, \op{deg}f_{2}=k$. Note that if $f_{1}(x), f_{2}(x)$ is a solution of (\ref{k22eq}), then due to the identity $\binom{x}{2}=\binom{1-x}{2}$, $f_{1}(1-x), f_{2}(1-x)$ is also a solution. In the sequel we count such pairs of solutions as one. We are motivated by findings presented in \cite{BBW}.

\begin{thm}\label{thmsol}
	Let $x$ be a variable.
	\begin{enumerate}
		\item[(1)] For $k=3, 5$ equation {\rm (\ref{k22eq})} has exactly three solutions.
		\item[(2)] For $k=7$ equation {\rm (\ref{k22eq})} has exactly one solution.
		\item[(3)] For $k\in\{9, 11, 13, 15, 17, 19\}$ equation {\rm (\ref{k22eq})} has no solutions.
	\end{enumerate}
\end{thm}

\section{Proofs of the theorems}
\begin{proof}[Proof of Theorem \ref{kl24}]
		In order to get the result it is enough to note that the equation $\binom{y}{2}=\binom{x}{4}+d$ can be rewritten as
		$$
		X^2-3Y^2=-2(12d+1),
		$$
		where $X=x^2-3x+1, Y=2y-1$. If $2(12d+1)\equiv 0\pmod{p}$, then $X^2\equiv 3Y^2\pmod{p}$. Under our assumption on $p$ we see that 3 is quadratic non-residue modulo $p$ and congruence \eqref{cong}, and hence equation \eqref{binomd}, has no integer solutions.
\end{proof}
	
	Motivated by the result above, we performed numerical search for pairs $(k,l), k\leq l\leq 10, d\in\Z$ and prime numbers $p>l$ such that the congruence (\ref{cong}) has no solutions modulo $p$. Here are results of our computations.
	{\tiny\begin{equation*}
	\begin{array}{|l|l|l||l|l|l|}
	\hline
	(k,l)  & p & d\pmod{p}   & (k,l)  & p  & d\pmod{p}              \\
	\hline
	(2,6)  & 7  & 4           & (6,8)  & 11 & 4                       \\
	(2,8)  & 11 & 7           &        & 13 & 3                       \\
	& 13 & 11          &        & 19 & 4                       \\
	(2,9)  & 11 & 8           & (6,9)  & 11 & 3, 4, 9                 \\
	(2,10) & 11 & 7, 8        & (6,10) & 11 & 2, 3, 4                 \\
	& 13 & 11          &        & 13 & 10, 11                  \\
	(3,4)  & 5  & 2           &        & 19 & 2, 4                    \\
	(3,8)  & 11 & 5           & (7,8)  & 11 & 4, 6                    \\
	(3,10) & 11 & 5           &        & 17 & 11                      \\
	(4,4)  & 5  & 2, 3        &        & 19 & 15                      \\
	(4,5)  & 7  & 3           & (7,9)  & 11 & 5, 6                    \\
	(4,6)  & 7  & 2, 3        & (7,10) & 11 & 4, 5, 6                 \\
	& 13 & 10          &        & 13 & 6                       \\
	& 19 & 2           & (8,8)  & 11 & 4, 5, 6, 7              \\
	(4,8)  & 11 & 8, 9        &        & 13 & 2, 11                   \\
	& 13 & 10, 11      &        & 17 & 4, 13                   \\
	(4,9)  & 11 & 7, 8        &        & 19 & 3, 16                   \\
	& 13 & 7           &        & 23 & 7, 16                   \\
	(4,10) & 11 & 6, 7, 8, 9  & (8,9)  & 11 & 3, 4, 5, 7              \\
	& 13 & 6, 10       & (8,10) & 11 & 2, 3, 4, 5, 6           \\
	& 23 & 9           &        & 13 & 4, 7, 10                \\
	(5,5)  & 7  & 3, 4        &        & 19 & 16                      \\
	& 11 & 3, 8        & (9,10) & 11 & 2, 3, 4, 5, 6, 7, 8     \\
	(5,6)  & 7  & 2, 3, 4     &        & 13 & 4, 6, 7, 8              \\
	& 11 & 2, 7, 8     &        & 17 & 8, 11, 14               \\
	(5,8)  & 11 & 5           &        & 23 & 7                       \\
	(5,9)  & 11 & 3, 8        & (10,10)& 11 & 2, 3, 4, 5, 6, 7, 8     \\
	(5,10) & 11 & 2, 3, 7, 8  &        & 13 & 4, 5, 6, 7, 8, 9        \\
	(6,6)  & 7  & 2, 3, 4, 5  &        & 17 & 5, 8, 9, 12             \\
	& 11 & 2, 3, 8, 9  &        & 19 & 3, 5, 14, 16            \\
	& 13 & 3, 10       &        & 23 & 5, 18                   \\
	&  &               &        & 29 & 6, 23                   \\ \hline
	\end{array}
	\end{equation*}}
\begin{center}
Table 2. Pairs $(k,l), k\leq l\leq 10$ such that there exist $p\in\mathbb{P},p\geq\op{max}\{k,l\}$ such that for some $d\in\{1,\ldots,p\}$ the congruence (\ref{cong}) has no solutions.
\end{center}
\begin{proof}[Proof of Theorem \ref{kl_eq}]
Here we obtain that
$$
\prod_{i=0}^{k-1}(n-i)-\prod_{i=0}^{k-1}(m-i)=d\cdot k!,
$$
and the polynomial is reducible. It follows that
$$
(n-m)F(n,m)=d\cdot k!.
$$
Hence $(n-m)$ divides $d\cdot k!.$ It remains to solve the one variable polynomial equation
$$
F(m+d_1,m)-\frac{d\cdot k!}{d_1}
$$
for $d_1|(d\cdot k!).$
\end{proof}

	\begin{rem} Let us note that if $k=l>2$, then in the considered range, i.e, $d\in\{-20,\ldots,20\}$ we have found at most one integer solution. It is an interesting problem to look for values of $d$ such that the equation
		\begin{equation}\label{keql}
		\binom{n}{k}=\binom{m}{k}+d
		\end{equation}
		has more than one solution in positive integers $m, n$ satisfying $n>m$.
		In order to construct values of $d$ such that equation (\ref{keql}) has ``many'' solutions we used the following strategy. First, we computed the set
		$$
		D_{k}:=\left\{\binom{n}{k}-\binom{m}{k}:\;k<m<n\leq 10^{4}\right\},
		$$
		and then looked for duplications in $D_{k}$. We considered $k\in\{3,\ldots, 10\}$. As one could expect, in the case $k=3$ the number of duplicates is big. In fact, we found 488 values of $d$ which appeared at least three times in $D_{3}$. The smallest value correspond to $d=2180$ with the solutions $(n,m)=(25,10), (33,28), (36,32)$. We found only three values of $d$ such that equation (\ref{keql}) has four solutions. The values of $d$ and the corresponding solutions are as follows:
		\begin{equation*}
		\begin{array}{lll}
		d=10053736   &  & (n,m)=(398,132), (628,572), (968,946), (990,969), \\
		d=209920964  &  & (n,m)=(1081,58), (1144,617), (1242,868), (3532,3498), \\
		d=1928818640 &  & (n,m)=(2266,362), (2268,428), (3622,3300), (4991,4831) .
		\end{array}
		\end{equation*}
		We strongly believe that the following is true.
		
		\begin{conj}
			For each $N\in\N$ there is $d_{N}\in\N$ such that the equation $\binom{n}{3}-\binom{m}{3}=d_{N}$ has at least $N$ positive integer solutions.
		\end{conj}
		
		For $k=4$ we found 1190 values of $d$ which appeared at least two times in $D_{4}$. The smallest value corresponds to $d=680$ with the solutions $(n,m)=(13,7), (18,17)$. We found only one value of $d$ such that equation (\ref{keql}) has three solutions. More precisely, for $d=18896570$ equation (\ref{keql}) has three solutions $(n,m)=(185,163), (258,251), (486,485)$.
		
		For $k=5$ we found 4 values of $d$ which appeared at least 2 times in $D_{5}$. The values of $d$ and the corresponding solutions are as follows:
		\begin{equation*}
		\begin{array}{lll}
		d=146438643  &  & (n,m)=(117,78), (133,118), \\
		d=153852348  &  & (n,m)=(118,78), (133,117), \\
		d=817514347  &  & (n,m)=(160,53), (209,197),\\
		d=2346409884 &  & (n,m)=(197,53), (209,160).
		\end{array}
		\end{equation*}
		
		For $k=6$ we also found 4 values of $d$ which appeared at least 2 times in $D_{6}$. The values of $d$ and the corresponding solutions are as follows:
		\begin{equation*}
		\begin{array}{lll}
		d=3819816      &  & (n,m)=(40,18), (57,56), \\
		d=32449872     &  & (n,m)=(56,18), (57,40), \\
		d=66273157776  &  & (n,m)=(193,66), (252,243),\\
		d=268624373556 &  & (n,m)=(243,66), (252,193).
		\end{array}
		\end{equation*}
		
		For $k=7$ we found only one value of $d\in D_{7}$ such that equation $(\ref{keql})$ has two solutions. For $d=8008$ we have solutions $(n,m)=(16,14), (17,16)$.
		
		For $k=8, 9, 10$ there are no duplicates in the set $D_{k}$.
	\end{rem}

\begin{proof}[Proof of Theorem \ref{elliptic}]
All the equations related to this part can be reduced to elliptic curves given is some model.
\begin{center}
	\tiny
	\begin{tabular}{|c|c|c|}
		\hline
		$(k,l)$	& equation & transformation \\
		\hline
		(2,3) & $Y^2=X^3 - 36X^2 + 288X + 10368d+1296$ &  $X=12m, Y=216n - 108$\\
		\hline
		(2,4)& $Y^2=3X(X-1)(X-2)(X-3)+72d+9$ & $X=m, Y=6n-3$ \\
		\hline
		(2,6)& $Y^2=X(X+40)(X+60)+10^4\cdot (72d+9)$ & $X=10m^2-50m, Y=600n-300$ \\
		\hline
		(2,8)& $Y^2=35X(X+6)(X+10)(X+12)+420^2(8d+1)$ & $X=m^2-7m, Y=420(2n-1)$ \\
		\hline
		(3,4)& $Y^2=X(X-4)(X-8)-384d+16$ & $X=4n, Y=4m^2-12m+4$ \\
		\hline
		(3,6)& $15X(X-1)(X+1)=Y(Y-3)(Y+4)+90d$  & $X=n-1, Y=(m-2)(m-3)/2$ \\
		\hline
		(4,6)& $Y^2=X(X+120)(X+180)+30^4\cdot (24d+1)$ & $X=30m^2-150m, Y=900(n^2-3n+1)$ \\
		\hline
		(4,8)& $Y^2=105X(X+6)(X+10)(X+12)+420^2(24d+1)$ & $X=m^2-7m, Y=420(n^2-3n+1)$ \\
		\hline
	\end{tabular} 	
\end{center}
\begin{center}
Table 3. Elliptic models of certain Diophantine equations of the form $\binom{m}{k}=\binom{n}{l}+d$
\end{center}
There exists a number of software implementations for finding integral points on elliptic curves~\cite{MAGMA,SAGE}.
These procedures are based on a method developed by Stroeker and Tzanakis \cite{StTz} and independently by Gebel, Peth\H{o} and Zimmer \cite{GPZell}. One may follow the transformations provided in \cite{StroekerWeger} to handle these cases. Here we used the Magma procedures \texttt{IntegralPoints()} and \texttt{IntegralQuarticPoints()}. In some cases there exist no solution and we used  \texttt{IsLocallySolvable()} and \texttt{TwoCoverDescent()} \cite{BStwo}. In cases related to $(k,l)=(3,6)$ we follow the above mentioned elliptic logarithm method, the cases with $d=-1,0,1$ were solved earlier as given in the introduction, so it remains to deal with the values $d\in\{-3,-2,2,3\}.$

The case $d=2$ yields an elliptic curve with Mordell-Weil rank 3 while the remaining
three values of $d$ yield elliptic curves with Mordell-Weil rank 2; we only provide
details for the case $d=2$.

For this case we set $u=X$, $Y=v$ and we have the equation
\begin{align}\label{g_d=2}
C: g(u,v)=0, \quad \mbox{where}\quad g(u,v)=15u^3-v^3+4v^2-15u-3v-180,
\end{align}
where $u=n-1$ and $v=\dfrac{1}{2}\left(\left(m-\dfrac{5}{2}\right)^2-\dfrac{1}{4}\right)=
 (m-2)(m-3)/2$
and the Weierstrass model which is birationally equivalent to $C$ over $\Q$ is
\begin{align}\label{fd2}
E: y^2=x^3-1575x-48749850=: f(x).
\end{align}

A notation remark: We will use ``exponents'' $^C$ and $^E$ on a point to declare
whether the point is viewed as one on $C$ or $E$, respectively. Also, we will use
$(u,v)$ or $(x,y)$ for the $C$-coordinates or the $E$-coordinates, respectively.

As already mentioned, $E(\Q)$ has rank 3; its free part is generated by the points
\[
 P_{1}^{E}=\left(10905/4, -1137285/8\right),\;
 P_{2}^{E}=\left(7465/9, 616040/27\right),\;
P_{3}^{E}=\left(10246/25, -551206/125\right)
\]
and the torsion subgroup is trivial.

The birational transformation between the models $C$ and $E$ is
\begin{align*}
C\ni P^C:=(u,v)&\longrightarrow (x,y)=\left(\mathcal{X}(u,v),
\mathcal{Y}(u,v)\right):=P^E\in E \\
C\ni P^C=:\left(\mathcal{U}(x,y),\mathcal{V}(x,y)\right)=(u,v)
 &\longleftarrow(x,y)=:P^E\in E
\end{align*}
with
\[
\mathcal{X}(u,v)= \frac{3(620u^2+235uv+106v^2-210u-438v+1960)}{(u+4)^2},
\]
\[
\mathcal{Y}(u,v)=
\frac{\mathcal{Y}\mbox{num}(u,v)}
{(u+4)^3},
\]
where
\begin{align*}
\mathcal{Y}\mbox{num}(u,v)=&3(45795u^3+19080u^2v+7285uv^2-35895u^2-16795uv-4568v^2+\\
&+32940u+65744v-408000)
\end{align*}
and
\begin{eqnarray}
\mathcal{U}(x,y) & = & \frac{4x^3-465x^2+318xy+3903030x-94455y+257567175}
 {-x^3+5580x^2-290250x+161614575}, \nonumber
 \\ & & \label{UVd2} \\
\mathcal{V}(x,y) & = & \frac{9x^3+7020x^2+705xy-9215775x+205560y+1359589050}
   {-x^3+5580x^2-290250x+161614575}. \nonumber
\end{eqnarray}

With the aid of Maple we find out that there is exactly one conjugacy class of Puiseux series $v(u)$ solving $g(u,v)=0$.
This unique class contains exactly three series and only the following
one has real coefficients:

\begin{align}\label{Puiseux_d_2}
v_1(u)= &\nonumber\,\zeta u+\frac{4}{3}+\left(\dfrac{7}{135}\zeta^2
-\dfrac{1}{3}\zeta\right)u^{-1}+\dfrac{968}{443}\zeta u^{-2}
+\left(\dfrac{7}{405}\zeta^2-\dfrac{1}{9}\zeta\right)u^{-3}
\\  &+\left(\dfrac{6776}{32805}\zeta^2-\dfrac{1936}{729}\zeta\right)u^{-4}+\ldots.
\end{align}

Here $\zeta$ is the cubic root of $15$. For every real solution of $g(u,v)=0$ with $|u|\geq 3$ it is true that
$v=v_1(u)$ (according to Lemma 8.3.1 in \cite{Tzanakis}).

Then the point $P_0^E$ that plays a crucial role in the resolution
(see \cite[ Definition 8.3.3]{Tzanakis}) is
\[
P_{0}^{E}=(318\zeta^2+705\zeta+1860, 21855\zeta^2+57240\zeta+137385).
\]

\noindent
Referring to the discussion of Section 1 of \cite{Katsipis}, we consider
the linear form
\label{page def L(P) when d=2}
\[
L(P)=\left(m_0+\dfrac{s}{t}\right)
\omega_1+m_1\ellog(P_1)+m_2\ellog(P_2)
+m_3\ellog(P_3)
\pm\ellog(P_0).
\]

Since $f(X)$ has only one real root, namely $e_1\approx366.7439448002$,
we have $E(\R)=E_0(\R)$, therefore $\mathfrak{l}(P_i)$ coincides with the
elliptic logarithm of $P_i^E$ for $i=1,\ldots,3$
(see Chapter 3 of \cite{Tzanakis}, especially,
Theorem 3.5.2).
On the other hand, $P_0^E$ has irrational coordinates.
As Magma does not possess a routine for calculating elliptic
logarithms of non-rational points, we wrote our own routine in Maple
for computing $\ellog$-values of points with algebraic coordinates.
\label{page my maple routine}
Thus we compute
\[
\mathfrak{l}(P_1)\approx0.0191558345, \quad
\mathfrak{l}(P_2)\approx-0.0349501519,\]
\[\mathfrak{l}(P_3)\approx0.0532999952,\quad
\mathfrak{l}(P_0)\approx-0.00763363355.
\]

Note that the four points $P_i^E,\,i=0,1,\ldots,3$  are $\Z$-linearly
independent because their regulator is non-zero
(see \cite[Theorem 8.1]{SchZim}).
Therefore our linear form $L(P)$ falls under the scope of the
second ``bullet'' in \cite[page 99]{Tzanakis}
and we have $r_0=1$, $s/t=s_0/t_0=0/1=0$, $d=1$, $r=4$, $n_i=m_i$ for
$i=1,\ldots,3$, $n_4=\pm1$, $n_0=m_0$, $k=r+1=4$,
$\eta=1$ and $N=\max_{0\leq i\leq 4}|n_i|
\leq r_0\max\{M,\frac{1}{2}rM+1\}+\frac{1}{2}\eta r_0
=\frac{3}{2}M+\frac{3}{2}$,
so that, in the relation (9.6) of \cite{Tzanakis} we can take
\begin{equation}
                \label{eq d=2, alpha,beta}
\alpha=3/2, \beta=3/2.
\end{equation}

We compute the canonical heights of $P_1^E, P_2^E, P_3^E$
using Magma\footnote{For the definition of the
canonical height we follow J.H.~Silverman; as a consequence the values
displayed here for the canonical heights are the halves of those computed
by Magma and the least eigenvalue $\rho$ of the height-pairing matrix
$\mathcal{H}$ below, is half that computed by Magma;
cf.~``Warning'' at bottom of p.~106 in \cite{Tzanakis}.
\label{foot canonical height}}
and for the canonical height of $P_0^E$ we confine ourselves to the upper
bound by applying \cite[Proposition 2.6.4]{Tzanakis}. Thus we have
\[
\hat{h}(P_1^E)\approx3.6037959076, \quad
\hat{h}(P_2^E)\approx3.7072405585,\]
\[\hat{h}(P_3^E)\approx4.8663287093,\quad
\hat{h}(P_0^E)\leq 8.022765298\,.\]

The corresponding height-pairing matrix for the particular Mordell-Weil
basis is
\[
\mathcal{H}\approx\left(\begin{array}{rrr}
3.6037959076 & -1.0424191872 & -1.2722619781 \\
-1.0424191872 & 3.7072405585 & 3.0174040388 \\
-1.2722619781 & 3.0174040388 & 4.8663287093
\end{array}\right)
\]
with minimum eigenvalue
\begin{equation}
       \label{eq d=2,rho}
  \rho\approx1.2142056695.
\end{equation}
Next we apply \cite[Proposition 2.6.3]{Tzanakis} in order to compute
a positive constant $\gamma$
with the property that
$\hat{h}(P^E)-\frac{1}{2}h(x(P))\leq \gamma $ for every point
$P^E=(x(P),y(P))\in E(\Q)$, where $h$ denotes Weil height;\footnote{In the notation of
\cite[Proposition 2.6.3]{Tzanakis},
as a curve $D$ we take the minimal model of $E$ which is $E$ itself.}
it turns out that
\begin{equation}
    \label{eq d=2, gamma}
\gamma\approx 4.8726444820.
\end{equation}

Finally, we have to specify the constants $c_{12},c_{13},c_{14},c_{15}$
defined in \cite[Theorem 9.1.2]{Tzanakis}.
This can be carried out almost automatically
with a Maple program. In this way we compute
\begin{equation}
                \label{eq d=2, c12,c13,c14,c15}
c_{12}\approx1.07690\cdot10^{27},\quad
c_{13}\approx4.04702\cdot10^{162},\quad\\
c_{14}\approx2.09861,\quad c_{15}\approx24.99686.
\end{equation}

 According to \cite[Theorem 9.1.3]{Tzanakis}, applied
to ``case of Theorem 8.7.2'', if $|u(P)|\geq\max\{B_2,B_3\}$, where
$B_2$ and $B_3$ are explicit positive constants, then either
$M\leq c_{12}$, where $c_{12}$ is an explicit constant, or
\begin{equation}
                \label{eq furnishes ub for M}
\rho M^2\leq \frac{c_{11}c_{13}}{2\theta}
(\log(\alpha M+\beta)+c_{14})
(\log\log(\alpha M+\beta)+c_{15})^{r+3}+\gamma
+\frac{c_{11}}{2\theta}\log\frac{c_9}{1+\theta}+
\textstyle{\frac{1}{2}}c_{10},
\end{equation}
where all constants involved in it are explicit.
More specifically, (in a similar way as
in Appendix B in \cite{Katsipis} for the case of $d=(N^3-N)/6$),
\[
B_2=4,\quad B_3=5,\quad
\theta=1,\quad c_9=0.17,\quad c_{10}=\log(11800),\quad c_{11}=2.
\]

So, in view of  \eqref{eq furnishes ub for M} and
\eqref{eq d=2, alpha,beta}, \eqref{eq d=2,rho},
\eqref{eq d=2, gamma}, \eqref{eq d=2, c12,c13,c14,c15},
we conclude that, if $|u(P)|\geq 5$, then either $M\leq c_{12}$ or
\begin{eqnarray*}
\lefteqn{1.2142056695\cdot M^2 \leq}
\\
 & & 4.04\cdot10^{162}\cdot(\log(1.5M+1.5)
+2.0986)\cdot(\log(\log(1.5M+1.5))+24.9968)^6+7.09542.
\end{eqnarray*}

But for all $M\geq 6.64\cdot10^{86}$, we check that the left-hand side
is strictly larger than the right-hand side which implies that
$M < 6.64\cdot10^{86}$, therefore
\begin{equation}\label{BM_d_2}
M\leq\max\{c_{12},\; 6.64\cdot10^{86}\}=6.64\cdot10^{86}
\quad\mbox{provided that $|u(P)|\geq 5$.}
\end{equation}

An easy straightforward computation shows that $P^C=(-4,-9)$
is the only one integer point with $|u(P)|\leq 4$ (equivalently, the integer solution $(u,v)$ of
\eqref{g_d=2} with $|u|\leq 4$).

In order to find explicitly all points $P^C$ with $|u(P)|\geq 5$ it is
necessary to reduce the huge upper bound \eqref{BM_d_2} to an upper
bound of manageable size. This is accomplished  with LLL-algorithm \cite{intLLL},
in a similar way as in Appendix D in \cite{Katsipis}, and we obtain the
reduced bound $M\leq 10$.
Therefore, we have to check which points
\[
P^E=m_1P_1^E+m_2P_2^E+m_3P_3^E,
   \quad\mbox{with $\max_{1\leq i\leq 3}|m_i|\leq 10,$}
\]
have the property that $P^E=(x,y)$ maps via the transformation
\eqref{UVd2} to a point $P^C=(u,v)\in C$ with integer
coordinates.
We remark here that every point $P^C$ with $u(P)$ integer
and $|u(P)|\geq 5$ is obtained in this way, but the converse is not
necessarily true; i.e. if $\max_{1\leq i\leq 3}|m_i|\leq 10$ and
the above $P^E$ maps to $P^C$ with integer coordinates, it is not
necessarily true that $|u(P)|\geq 5$.
\noindent
 After a computational search we find the only one point $P^C=(-4,-9)$ which corresponds to the zero point $\mathcal{O}\in E$.

 So no integral solution $(m,n)$ (with $n\geq k$ and $m\geq l$) of equation \eqref{binomd} with $d=2$ exists.

For the other three cases we provide some details in the tables below:

\begin{equation*}
\begin{array}{|l|l|l|l|l|l|}
\hline
  d & a(d) & r  & \mbox{Generators}  & \rho & e_1 \\
  \hline
  -2 & -49559850 & 2 & P_1=(956289/4, 935155287/8)      & 1.499191 & 368.748212 \\
     &           &   & P_2=(198006/169, -86688954/2197) &          &            \\
  -3 & -111271725 & 2 & P_1=(1230, -41805)             & 2.568215 & 482.072907 \\
     &            &   & P_2=(221597697975/91145209,    &           & \\     
     &            &   & \quad\quad  103896688780607535/870163310323) &   &  \\                       
   3 & -110056725 & 2 & P_1=(1072825/2116, -429530005/97336) & 1.786872 & 480.319851 \\
     &            &   & P_2=(16866855/34969,                 & &\\
     &            &   & \quad\quad -7734674565/6539203)     & & \\    
  \hline
\end{array}
\end{equation*}
\begin{center}
Table 4. $C: 15u^3-v^3+4v^2-15u-3v-90d$ and $E: y^2=x^3-1575x+a(d)$
\end{center}

\begin{equation*}
\begin{array}{|c|c|c|}
\hline 
d  & B(M): \mbox{Initial bound} & \mbox{Reduced bound} \\
\hline
  -2  & 5.06\cdot10^{62} & 6  \\ 
  -3  & 8.66\cdot10^{62} & 5 \\
  3   & 9.07\cdot10^{62}  & 5\\
\hline
\end{array}
\end{equation*}
\begin{center}
Table 5. Upper bounds of $M$.
\end{center}

\begin{equation*}
\begin{array}{|c|c|c|}
\hline
d & P^E=(x,y) & P^C=(u,v) \\
\hline
  -2  & \mathcal{O} & (-2, 6) \\
  -3  & \mathcal{O} & (3, 10) \\
  3   & \mathcal{O},\; (16155, -2053305) & (3, 6),\;(4, 10)\\
\hline
\end{array}
\end{equation*}
\begin{center}
Table 6. All points $P^E=\Sigma_i m_iP_i^E$ with $P^C=(u,v)\in\Z\times\Z$.
\end{center}
\end{proof}

\begin{proof}[Proof of Theorem \ref{hyperell}]
We provide details only in case of $d=3,$ here the rank of the Jacobian is 6 (like in case of $d=1$).
Equation \eqref{binomd} with $d=3$ defines the hyperelliptic curve
$$
y^2=15x(x-1)(x-2)(x-3)(x-4)+75^2.
$$
Based on Stoll's papers \cite{StollB1}, \cite{Stoll}, \cite{StollB2} one can determine generators for the Mordell-Weil group by using Magma \cite{MAGMA}. We obtain that $J(\mathbb{Q})$ is free of rank $6$ with
Mordell-Weil basis given by (in Mumford representation)
\begin{eqnarray*}
	&& D_1=<x - 4, -75>,\\
	&& D_2=<x - 3, 75>,\\
	&& D_3=<x - 1, -75>,\\
	&& D_4=<x, 75>,\\
	&& D_5=<x^2 - 7x + 30, 195>,\\
	&& D_6=<x^2 - 3x + 20, -30x - 45>
\end{eqnarray*}
and  the torsion subgroup is trivial. We apply Baker's method \cite{intBH} to get a large upper bound for $\log|x|,$ here we use the improvements given in \cite{BMSST} and \cite{GallHyp}. It follows that
$$
\log|x|\leq 1.028\times 10^{612}.
$$
We have from Corollary 3.2 of \cite{GallHyp} that every integral point on the curve can be
expressed in the form
$$
P-\infty=\sum_{i=1}^{6}n_iD_i
$$
with $||(n_1,n_2,n_3,n_4,n_5,n_6)||\leq 1.92\times 10^{306}=:N.$ Proposition 6.2 in \cite{GallHyp} gives an estimate for the precision we need to compute the appropriate matrices, this bound is as follows
$$
((1/5)(48\sqrt{r}Nt + 12\sqrt{r}N + 5N + 48))^{(r+4)/4}\approx 2.6\times 10^{769},
$$
where in our case $r=6$ and $t=1.$ We choose to compute the period matrix and the  hyperelliptic  logarithms with 1500 digits of precision. The hyperelliptic logarithms of the divisors $D_i$ are given by
\begin{eqnarray*}
\varphi(D_1)&=&(0.087945\ldots+i0.112834\ldots,-0.473844\ldots-i0.741784\ldots)\in\mathbb{C}^2,\\
\varphi(D_2)&=&(0.114612\ldots+i0.112834\ldots,-0.420527\ldots-i0.741784\ldots)\in\mathbb{C}^2,\\
\varphi(D_3)&=&(-0.044486\ldots+i1.333456\ldots,-0.416321\ldots+i5.329970\ldots)\in\mathbb{C}^2,\\
\varphi(D_4)&=&(0.127905\ldots+i0.112834\ldots,-0.413878\ldots-i0.741784\ldots)\in\mathbb{C}^2,\\
\varphi(D_5)&=&(-0.118415\ldots+i0.037611\ldots,-0.857076\ldots-i0.247261\ldots)\in\mathbb{C}^2,\\
\varphi(D_6)&=&(0.128537\ldots+i0.075223\ldots,-0.173077\ldots-i0.494522\ldots)\in\mathbb{C}^2.
\end{eqnarray*}
We need now to choose an integer $K$ that is larger than the constant given by Proposition 6.2 in \cite{GallHyp}. Setting $K=10^{1300}$ we get a new bound $126.98$ for $||(n_1,n_2,n_3,n_4,n_5,n_6)||.$ We repeat the reduction process with $K=10^{16}$ that yields a better bound, namely $15.6.$ Two more steps with $K=6\times 10^{11}$ and $K=2\times 10^{11}$ provide the bounds $13.94$ and $13.8.$ It remains to compute all possible expressions of the form
$$n_1D_1+\ldots+n_6D_6$$
with $||(n_1,n_2,n_3,n_4,n_5,n_6)||\leq 13.8.$ We performed a parallel computation to enumerate linear combinations coming from integral points on a machine having 12 cores. The computation took 3 hours and 23 minutes. We obtained the following non-trivial solutions
\begin{eqnarray*}
\binom{11}{5}+3&=&\binom{31}{2},\\
\binom{16}{5}+3&=&\binom{94}{2},\\
\binom{375}{5}+3&=&\binom{346888}{2},\\
\binom{379}{5}+3&=&\binom{356263}{2}.
\end{eqnarray*}

If $d=1,$ then the rank of the Jacobian is 6 and the Baker bound is $\log |x|\leq 1.225\times 10^{532}$ and we have that $||(n_1,n_2,n_3,n_4,n_5,n_6)||\leq 2.23\times 10^{266}.$ In three steps it is reduced to $14.97.$ In this case the non-trivial solutions are as follows
\begin{eqnarray*}
	\binom{10}{5}+1&=&\binom{23}{2},\\
	\binom{22}{5}+1&=&\binom{230}{2},\\
	\binom{62}{5}+1&=&\binom{3598}{2},\\
	\binom{135}{5}+1&=&\binom{26333}{2},\\
	\binom{139}{5}+1&=&\binom{28358}{2}.
\end{eqnarray*}

If $d=-3,-1,2,$ then the rank of the Jacobian is 3, we followed the arguments given in \cite{BMSST} and \cite{Gallegos} to obtain a large bound for the size of possible integral solutions. We present them in the table below.
\begin{equation*}
	\begin{array}{|c|c|}
		\hline
		d & \mbox{bound for}\;\log |x|\\
		\hline
		-3 & 2.91\cdot 10^{608}\\
		-1 & 1.21\cdot 10^{552}\\
		2 & 3.25\cdot 10^{590} \\
		\hline
	\end{array}
\end{equation*}
\begin{center}
Table 7. Upper bounds for $\log|x|$.
\end{center}

In all three cases the rank of the Jacobians are equal to 3 and the torsion subgroup is trivial hence all points can be written as
$$
n_1D_1+n_2D_2+n_3D_3,
$$
where $n_i\in\mathbb{Z}.$ Using the previously applied hyperelliptic logarithm method the initial large upper bounds for $\max\{|n_i|\}$ can be significantly reduced. If $d=-3,$ then after one reduction step we get the bound 64 and other two steps make it 7. The only pair of integral points we obtain is given by $(6,\pm 75).$ Therefore we have
$$
\binom{3}{2}=\binom{6}{5}-3.
$$
If $d=-1,$ then first we obtain a reduced bound 51 and finally it follows that $\max\{|n_i|\}\leq 5.$ The complete list of integral points is given by $(5,\pm 15),(8,\pm 315).$ Thus we obtain
$$
\binom{11}{2}=\binom{8}{5}-1.
$$
Finally, in case of $d=2$ the first reduction yields a bound 58 and the third one provides 6. The complete set of integral solutions is $\{(-1,\pm 45),(5,\pm 75)\},$ so we do not get non-trivial solution of \eqref{binomd}.

If $d=-2,$ then the rank of the Jacobian is 1, therefore classical Chabauty's method \cite{Chab} can be applied, it is now implemented in Magma \cite{MAGMA}. We obtain that the equation $\binom{n}{2}=\binom{m}{5}-2$ has no non-trivial solution.
\end{proof}

	\begin{rem}
		Let
		$$
		C_{d}:\; y^2=15x(x-1)(x-2)(x-3)(x-4)+15^2(8d+1)
		$$
		and write $J_{d}:=\op{Jac}(C_{d})$. The curve $C_{d}$ is isomorphic to the curve defined by the equation $\binom{y}{2}=\binom{x}{5}+d$.
		We computed upper bounds for the numbers $r_{d}=\op{rank}J_{d}(\Q)$ using the Magma procedure {\tt RankBound}.
		We obtained the following data
		\begin{equation*}
		\begin{array}{|l|l|}
		\hline
		i & \mbox{the value of $d$ such that}\; r_{d}\leq i\\
		\hline
		0 &-45,-40,-39,-37,-34,-10,-9,-4,8,25,26,40,47\\
		1 &-47,-36,-33,-31,-28,-26,-25,-22,-14,-13,-8,-5,-2,5,\\
		&11,17,20,29,32,41,50\\
		2 &-50,-46,-41,-38,-32,-30,-29,-24,-23,-19,-16,-7, 4, 13\\
		&14,23,30,31,38,43,44\\
		3 &-48,-44,-43,-42,-35,-21,-20,-15,-11,-3,-1,2,7,16,18\\
		&19,33,35,39,42,48\\
		4 &-49,-27,-18,-17,-12,-6,9,12,22,24,34,37,46,49\\
		5 &27,36\\
		6 &0,1,3,6,10,15,45\\
		7 &21,28 \\
		\hline
		\end{array}
		\end{equation*}
\begin{center}
Table 8. Upper bounds for the rank of Jacobian of the curve $C_{d}$ for $d\in\{-50,\ldots,50\}$.  
\end{center}
		We checked that for $i\in\{0,4,5,6,7\}$ the upper bounds computed by {\tt RankBound} are actually equal to the ranks. 
	
	Let us note that $21=\binom{7}{2}$ and $28=\binom{8}{2}$. We checked that in both cases the rank is equal to 7. This follows from the existence of seven independent divisors in $J_{d}(\Q)$. They are as follows:
	\begin{align*}
	d=21;\;&<x - 3, -345>, <x - 1, -345>, <x - 4, 345>, <x, 345>, \\
	&<x + 3, 285>, <x + 4, 135>, <x - 11, 975>, <x^2 + x + 30, -30x + 165>,\\
	d=28;\;&<x - 3, 225>, <x - 1, -225>, <x - 4, 225>, <x - 12, 1215>,\\
	& <x - 17, -3345>, <x, 225>, <x^2 - x + 18, -135>.
	\end{align*}
	We also looked for high rank Jacobians for further values of $d$ of the form $\binom{w}{2}.$
	For $d=66=\binom{12}{2}$ we obtained the equality $r_{66}=8$ with the following independent divisors
	\begin{align*}
	&<x - 3, -345>, <x - 1, -345>, <x - 4, 345>, <x, 345>, \\
	&<x + 3, 285>, <x + 4, 135>, <x - 11, 975>, <x^2 + x + 30, -30x + 165>.
	\end{align*}
	The torsion part of $J_{66}(\Q)$ is trivial. We conjecture that the only solutions in positive integers of the equation $\binom{y}{2}=\binom{x}{5}+66$ are
	\begin{align*}
	(x,y)=&(1,23), (2,23), (3,23), (4,23), (11,65), (28,887), \\
	&(7935,1447264765), (7939,1449089815).
	\end{align*}
	The large points are explained by the fact that on the curve $C_{\binom{w}{2}}$ we have the following solutions
	\begin{eqnarray*}
	x&=&3\cdot5\cdot(2w-1)^2,\\
	y&=&75(720w^4 - 1440w^3 + 1020w^2 - 300w + 31)(2w-1) \mbox{ and }\\
	x&=&3\cdot5\cdot(2w-1)^2+4,\\
	y&=&75(720w^4 - 1440w^3 + 1140w^2 - 420w + 61)(2w-1).
    \end{eqnarray*}
Hence we obtain the following divisors on $J_{\binom{w}{2}}(\Q)$
\begin{eqnarray*}
&& (x, 30w - 15, 1),\\
&& (x - 1, 30w - 15, 1),\\
&& (x - 2, 30w - 15, 1),\\
&& (x - 3, 30w - 15, 1),\\
&& (x - 4, 30w - 15, 1),\\
&& (x - 60w^2 + 60w - 15, 108000w^5 - 270000w^4 + 261000w^3 - 121500w^2 + 27150w - 2325, 1),\\
&& (x - 60w^2 + 60w - 19, 108000w^5 - 270000w^4 + 279000w^3 - 148500w^2 + 40650w - 4575, 1).
\end{eqnarray*}
\end{rem}
	
\begin{rem}
In case of the equation
$$
\binom{n}{2}=\binom{m}{7}+d
$$
one obtains genus 3 curves. Stoll \cite{Stoll-g3} proved that the rank of the Jacobian is 9 if $d=0.$ For other values of $d$ in the range $\{-3,\ldots,3\}$ many of the genus 3 hyperelliptic curves have high ranks as well. Balakrishnan et. al. \cite{Bala-g3} developed an algorithm to deal with genus 3 hyperelliptic curves defined over $\mathbb{Q}$ whose Jacobians have Mordell-Weil rank 1. If $d=-2,$ then the equation is isomorphic to the curve
$$
Y^2=70X^7 - 1470X^6 + 12250X^5 - 51450X^4 + 113680X^3 - 123480X^2 + 50400X - 661500
$$
and using Magma (with \texttt{SetClassGroupBounds("GRH")} to speed up computation) we get that the rank of the Jacobian is 1. Therefore we may try to use the Sage implementation described in \cite{Bala-g3} to compute the set of rational points on this curve. The affine points are $(8,\pm 1470),$ hence we have the solution
$$
\binom{4}{2}=\binom{8}{7}-2.
$$
\end{rem}

\begin{proof}[Proof of Theorem \ref{thmsol}]
		In each case we will be working in the same way. More precisely, for given $k$ we write $f_{1}(x)=a_{2}x^2+a_{1}x+a_{0}$ and $f_{2}(x)=\sum_{i=0}^{k}b_{i}x^{i}$. The polynomial $\binom{f_{1}(x)}{k}+\binom{x}{2}-\binom{f_{2}(x)}{2}=\sum_{i=0}^{2k}A_{i}x^{i}$ needs to be zero. Thus the coefficient near $x^{i}$ in $F_{k}(x)$ need to be zero for $i=0,\ldots,2k$. In consequence, we are interested in solving the system of polynomial equations
		$$
		S_{k}:\;A_{0}=A_{1}=\ldots=A_{2k}=0
		$$
		in $k+4$ variables $a_{0}, a_{1}, a_{2}, b_{0},\ldots, b_{k}$. We have $A_{2k}=\frac{a_{2}^{k}}{k!}-\frac{b_{k}^{2}}{2}$ and thus $a_{2}=\frac{k!}{2}t^{2}, b_{k}=\left(\frac{k!}{2}\right)^{\frac{k-1}{2}}t^{k}$ for some non-zero $t\in\Q$. We note that after the substitution of the computed values of $a_{2}, b_{k}$ into the system $S_{k}$, the related system of equations
		$$
		S_{k}':\;A_{k}'=A_{k+1}'=\ldots=A_{2k-1}',
		$$
		where $A_{i}'$ comes from $A_{i}$ after the substitution of the computed values of $a_{2}, b_{k}$, is triangular with respect to the variables $b_{0}, b_{1}, \ldots, b_{k-1}$. More precisely, we have $\op{deg}_{b_{i}}A_{k+i}'=1$ for $i=0,\ldots, k-1$. Moreover, the coefficient near $b_{i}$ is a power of $t$ times a rational number. Solving for $b_{0},\ldots,b_{k-1}$ and substituting into $S_{k}'$ we are left with the system of equations
		$$
		S_{k}'':\;A_{0}''=A_{1}''=\ldots=A_{k-1}'',
		$$
		in three variables $a_{0}, a_{1}, t$. The polynomial $A_{i}''$ is the numerator of the rational function $A_{i}'$ after substitution of the computed values $b_{0},\ldots,b_{k-1}$. It seems that for each fixed odd $k\geq 3$, the system $S_{k}''$ can be solved using Gr\"{o}bner bases techniques. More precisely, we compute $G_{k}$ - the Gr\"{o}bner basis of the ideal generated by the polynomials $A_{i}'', i=0,\ldots,k-1$. For $k\geq 5$ we have more equations than variables we expect that the system $S_{k}''$ for all sufficiently large $k$ has no rational (and even complex) solutions. This can be confirmed with our approach for $k\in\{11,\ldots, 19\}$. However, we were unable to prove such a statement in full generality.
		
		We prove the first part of our theorem. However, we present details of the reasoning only for $k=3$. The case $k=5$ is proved in exactly the same way.  We are interested in rational solutions of the system
		$$
		S_{3}:\;A_{0}=\ldots=A_{6}=0.
		$$
		
		We have $a_{2}=3t^2, b_{3}=3t^3$ for some $t\neq 0$. We put the values of $a_{2}, b_{3}$ into the system $S_{3}$ and solve corresponding system of equations
$$
S_{3}':\;A_{3}'=A_{4}'=A_{5}'=0,
$$
 with respect to $b_{0}, b_{1}, b_{2}$. We get
		\begin{equation*}
		b_{0}=\frac{36 a_0 a_1 t^2-36 a_1 t^2-a_1^3+72 t^3}{144 t^3},\; b_{1}=\frac{12 a_0 t^2+a_1^2-12 t^2}{8 t}, \;b_{0}=\frac{3 a_1 t}{2}.
		\end{equation*}
		In consequence, after the substitution of the values of $a_{2}, b_{0}, b_{1}, b_{2}, b_{3} $ into the system $S_{3}$ we obtain the system
		$$
		S_{3}'':\;A_{0}''=A_{1}''=A_{2}''=0,
		$$
		where $A_{i}''=t^{2(3-i)}A_{i}'\in\Q[t,a_{0},a_{1}]$. It is an easy task to solve the system $S_{3}''$. Indeed, we compute Gr\"{o}bner basis $G_{3}$, of the ideal generated by $A_{0}'', A_{1}'', A_{2}''$. The basis $G_{3}$ contains four polynomials. Two of them are the following
		$$
		a_1^5(a_1+3)(a_1+12),\;(4 a_0-7) a_1^5(a_1+12)
		$$
		and we easily obtain the following solutions
		\begin{equation*}
		\begin{array}{lll}
		f_{1}(x)=3 (-1 + 2 x)^2, &  & f_{2}(x)=2 - 15 x + 36 x^2 - 24 x^3, \\
		f_{1}(x)=5 - 12 x + 12 x^2, &  & f_{2}(x)=5 - 21 x + 36 x^2 - 24 x^3, \\
		f_{1}(x)=\frac{1}{4}(12 x^2-12 x+7), &  & f_{2}(x)=\frac{1}{8}(-24 x^3+36 x^2-18 x+7).
		\end{array}
		\end{equation*}
		Note that the first two solutions were presented in \cite{BBW}. Unfortunately, the polynomials from the third solution take only non-integer values.
		
		For $k=5$ we proceed in the same way and omit details. However, let us note that the Gr\"{o}bner basis $G_{5}$ contains 7 polynomials. Two of them are the following
		$$
		a_1^9(a_1+60)(3 a_1+80),\; a_{1}^9(3 a_0-26)(a_1+60)
		$$
		and we obtain two solutions with integer coefficients and the solution (corresponding to the triple $t=2/3, a_{0}=26/3, a_{1}=-80/3$)
		$$
		f_{1}(x)=\frac{2}{3}(40 x^2-40 x+13),\;f_{2}(x)=\frac{1}{27}(12800 x^5-32000 x^4+32000 x^3-16000 x^2+3955 x-364).
		$$
		By replacing $x$ by $3x-1$ we obtain polynomial with integer coefficients, which is exactly the third solution from the paper \cite{BBW}.
		
		For $k=7$ the Gr\"{o}bner basis $G_{7}$ contains 11 elements. In
		particular, the following three polynomials are in $G_{7}$:
		$$
		a_1^{12}(a_1+70),\;a_1^{12}(2 a_0-41),\;a_1^{10}(420 t-a_1)(a_1+420t).
		$$
		We found that the only solution (corresponding to $t=1/6, a_{0}=41/2, a_{1}=-70$) is the following
		\begin{align*}
		f_{1}(x)&=\frac{1}{2}(140 x^2-140 x+41),\\
		f_{2}(x)&=\frac{1}{96}(5488000 x^7-19208000 x^6+28812000 x^5-24010000 x^4+11997160 x^3-3589740 x^2+\\
		&+594370 x-41847).
		\end{align*}

		The last part of our theorem follows from certain Gr\"{o}bner basis computations. For $k\in\{9, 11, 13, 15, 17, 19\}$ we found that the $G_{k}$ contains  polynomial of the form $t^{u_{k}}$ for some $u_{k}\in\N_{+}$, i.e., $t$ need to be zero which leads to contradiction.
	\end{proof}
	
	\begin{rem}
		{\rm Using the same approach as in the proof of the above theorem one can prove that the Diophantine equation $\binom{f_{1}(x)}{k}-\binom{x}{2}=\binom{f_{2}(x)}{2}$ has no polynomial solutions $f_{1}, f_{2}\in\Q[x]$ satisfying $\op{deg}f_{1}=2, \op{deg}f_{2}=k$ for $k\in\{3, 5,\ldots, 19\}$.
			
			We also looked for solutions of the more general Diophantine equation
			\begin{equation}\label{k22eq2}
			\binom{f_{1}(x)}{k}+\binom{f_{0}(x)}{2}=\binom{f_{2}(x)}{2},
			\end{equation}
			where $f_{0}$ is of degree 2. By using the same approach as in the proof of Theorem \ref{thmsol} one can prove that for $k\in\{5,7,\ldots, 19\}$ there are no solutions $f_{0}, f_{1}, f_{2}\in\Q[x]$ of (\ref{k22eq2}) satisfying $\op{deg}f_{0}=\op{deg}f_{1}=2$ and $\op{deg}f_{2}=k$.
			
			However, if we allow $f_{0}$ to be of degree 3 we found the following solutions. For $k=5$ we have the solution
			\begin{align*}
			f_{1}(x)&=15x^2, \\
			f_{0}(x)&=\frac{1}{2} \left(30 x^3-5 x+1\right),\\
			f_{2}(x)&=\frac{1}{2} \left(225 x^5-75 x^3+7 x+1\right).
			\end{align*}
			
			For $k=7$ we have the solution
			\begin{align*}
			f_{1}(x)&=2520 x^2+1, \\
			f_{0}(x)&=\frac{1}{2}\left(17640 x^3-23 x+1\right),\\
			f_{2}(x)&=\frac{1}{2} \left(32006016000 x^7-88905600 x^5+52920 x^3+7 x+1\right).
			\end{align*}
			Note that in both cases by replacing $x$ by $2x-1$ we get polynomials with integer coefficients.
			
			Playing around with the Diophantine equation $\binom{f_{0}(x)}{3}+\binom{f_{1}(x)}{3}=\binom{f_{2}(x)}{2}$ we also found the polynomial solution
			$$
			f_{0}(x)=x(3x+2),\;f_{1}(x)=(2 x+1) (3 x+2), \;f_{2}(x)=9 x^3+15 x^2+6 x+1.
			$$
		}
	\end{rem}

\bibliography{allbib}
\end{document}